\documentclass[12pt]{article}       
 {
\setlength{\oddsidemargin}{0.in} {
\setlength{\textwidth}{6.in} {
\setlength{\textheight}{9.in} {

\usepackage{graphicx}

\usepackage{amsfonts}
\usepackage{amsmath}
\usepackage{amssymb}
\usepackage[applemac]{inputenc}

\usepackage{graphicx}
\usepackage{subfigure}
\usepackage{mathabx}

\newcommand{\vs}[1]{\vspace{#1}}

\newcommand{\begitem}{\begin{itemize}}
\newcommand{\finit}{\end{itemize}}
\newcommand{\begenum}{\begin{enumerate}}
\newcommand{\finenum}{\end{enumerate}}
\newcommand{\begar}{\begin{array}}
\newcommand{\finar}{\end{array}}
\newcommand{\begeq}[1]{\begin{equation} \label{#1}}
\newcommand{\fineq}{\end{equation}}
\newcommand{\begct}{\begin{center}}
\newcommand{\finct}{\end{center}}
\newcommand{\ra}{\rightarrow}

\newcommand{\zun}{\vs{0.1cm}} 

\newcommand{\zdeux}{\vs{0.2cm}} 
\newcommand{\ztrois}{\vs{0.3cm}}
\newcommand{\zcinq}{\vs{0.5cm}} 

\newcommand{\bN}{\mathbb{N}}
\newcommand{\bR}{\mathbb{R}}

\newcommand{\bP}{\mathbb{P}}

\newcommand{\bI}{\mathbb{I}}

\newcommand{\limprob}{ \stackrel{\bP}{\longrightarrow} }

\begin{document}

\newtheorem{theo}{Theorem}
\newtheorem{prop}{Proposition}
\newtheorem{defi}{Definition}
\newtheorem{lem}{Lemma}
\newtheorem{cor}{Corollary}
\newtheorem{rmk}{Remark}

\begin{center}\large
{\Large Moment estimators of the extreme value index for randomly censored data in the Weibull domain of attraction }
\vspace{0.5cm}

Julien Worms\footnote{\it
Universit\'e de Versailles-Saint-Quentin-en-Yvelines, Laboratoire de Math\'ematiques de Versailles (CNRS UMR 8100), F-78035 Versailles Cedex, France, {\scriptsize\tt julien.worms@uvsq.fr}
} , Rym Worms\footnote{\it
 Universit\'e Paris-Est, Laboratoire d'Analyse et de Math\'ematiques Appliqu\'ees (CNRS UMR 8050), UPEMLV, UPEC, F-94010, Cr\'eteil, France, {\scriptsize\tt rym.worms@u-pec.fr}
}
\vspace{0.5cm}

21 october 2014
\vspace{0.8cm}

\normalsize

{\it Abstract}
\zdeux

\parbox{12.cm}{This paper addresses the problem of estimating the extreme value index in presence of random censoring for distributions in the Weibull domain of attraction. The  methodologies introduced in \cite{WWExtremes2014},  in the heavy-tailed case, are adapted here to the negative  extreme value index framework, leading to the definition of weighted versions of the popular moments of relative excesses with arbitrary exponent . This leads to the definition of  two families of  estimators (with an adaptation of the so called Moment estimator as a particular case), for which the consistency is proved under a first order condition. Illustration of their performance, issued from an extensive simulation study, are provided.
}

\end{center}

\vspace{2.cm}

\zdeux

\newpage

\section{Introduction}
 \label{intro}

Extreme value statistics is an active domain of research, with numerous fields of application, and which benefits from an important litterature in the context of i.i.d. data, dependent data, and (more recently) multivariate or spatial data. By contrast, methodological articles in the case of randomly censored data are quite recent and few : \cite{Einmahl2008} presents a general method for adapting estimators of the extreme value index in a censorship framework (a methodology based on a previous work \cite{Beirlant2007}),  \cite{DupuyDiopNdao2014} extends the framework to data with covariate information, and \cite{WWExtremes2014} proposes a more survival analysis-oriented approach restricted to the heavy tail case. Other existing works on the topic of extremes for censored data are \cite{BrahimiMeraghniNecir2013} and the review paper \cite{GomesNeves2011}.
\ztrois

In this paper, the topic of extreme value statistics for randomly censored data with negative extreme value index is addressed.
Our initial purpose was to rely on the ideas of \cite{WWExtremes2014} in order to define a more "natural" version (with respect to that proposed in \cite{Einmahl2008})  of the moment estimator in the context of censored observations. We finally came out to propose weighted versions of the popular moments of the relative excesses (with arbitrary exponent), and therefore define competitive estimators of the extreme value index in this censoring situation, for distributions in the Weibull maximum domain of attraction. \zun

 Let us first define more precisely the framework, the data, and the notations.
\ztrois

In the classical univariate framework of i.i.d. data, a central task is to estimate the extreme value index $\gamma$, which captures the main information about the behavior of the tail distribution of the data. More precisely, a distribution function (d.f.) $F$ is said to be in the maximum domain of attraction of $H_{\gamma}$ (noted $F \in D(H_{\gamma})$) with
 \[
 H_{\gamma} (x) := \left\{ \begar[c]{ll}
  \exp \left(-(1+\gamma x)^{-{1/\gamma}}\right) & \mbox{ for  $\gamma \neq 0$ and  $1+\gamma x >0$} \\
  \exp(-\exp(-x)) &  \mbox{ for  $\gamma = 0$ and  $x \in \bR$ },  \finar\right.
 \]
if there exist two normalizing sequences $(a_n) \subset \bR^+$ and $(b_n) \subset \bR$ such that
$$
 F^n(a_n x +b_n) \ \stackrel{n\rightarrow\infty}{\longrightarrow} \ H_{\gamma}(x)   \hspace{0.5cm}
 (\forall x \in \bR).
$$

We consider in this paper two independent i.i.d. non-negative samples $(X_i)_{i \leq n}$ and  $(C_i)_{i \leq n}$ with respective continuous distribution functions $F$  and $G$ (with end-points $\tau_F$ and $\tau_G$, where $\tau_F:= \sup  \{x, \ F(x)<1    \} $). In the context of randomly right-censored observations, one  only observes, for $1 \leq i \leq n$,
$$
  Z_i = X_i \wedge C_i \mbox{ \ and \ } \delta_i= \bI_{X_i \leq C_i}.
$$
We denote by $H$ the distribution function of the $Z$-sample, satisfying
$$
 1-H=(1-F)(1-G)
$$
and by  $Z_{1,n} \leq \cdots \leq Z_{n,n}$  the associated order statistics. In the whole paper, $\delta_{1,n}, \ldots, \delta_{n,n}$  denote the $\delta$'s corresponding to $Z_{1,n}, \ldots, Z_{n,n}$, respectively.
 $F$ and $G$ are assumed to belong to the maximum domains of attraction  $D(H_{\gamma_X})$   and $D(H_{\gamma_C})$ respectively,
where $\gamma_X$ and $\gamma_C$ are real numbers, which  implies that $H  \in D(H_{\gamma})$, for some $\gamma \in \bR$.
\zdeux

Our goal is to estimate the extreme value index $\gamma_X$ in this context of right censorship. The most interesting cases, described in \cite{Einmahl2008}, are  the following  :
\begin{eqnarray*}
\mbox{ case 1: }  & \gamma_X >0 \ ,\  \gamma_C >0 & \mbox{ in this case  }   \gamma= \frac{\gamma_X \gamma_C}{\gamma_X+ \gamma_C} \label{cas1} \zun \\
\mbox{ case 2: }  & \gamma_X <0 \ ,\  \gamma_C < 0 \ , \  \tau_F=\tau_G    & \mbox{ in this case  }  \gamma= \frac{\gamma_X \gamma_C}{\gamma_X+ \gamma_C} \label{cas2} \zun \\
\mbox{ case 3: }  & \gamma_X = \gamma_C = 0 \ ,\  \tau_F=\tau_G=+\infty \hspace{0.4cm} &  \mbox{ in this case  }  \gamma= 0 . \label{cas3}
\medskip
\end{eqnarray*}

In \cite{WWExtremes2014},  case 1 above was considered and an adaptation of the so-called Hill estimator to the right censoring framework was proposed.  In this paper, our aim is to consider case 2 above and adapt the approach leading to the so-called Moment Estimator to this censored situation.  An adaptation of this estimator was already proposed in \cite{Einmahl2008} : it consists in dividing the  classical Moment Estimator $\hat{\gamma}^Z_n$  of $\gamma$ (calculated from the $Z$-sample) by the proportion $$ \textstyle \widehat{p}:=k_n^{-1}\sum_{i=1}^{k_n} \delta_{n-i+1,n}  $$ of uncensored data in the tail, where $k_n$ is the number of upper order statistics retained. Note that $\hat{\gamma}^Z_n$ is an appropriate  combination of the following moments
\[
{\cal M}^{(\alpha)}_{n,k_n} := \frac{1}{k_n} \sum_{i=1}^{k_n} \log^{\alpha}  \left(  \frac{Z_{n-i+1,n}}{Z_{n-k_n,n}} \right),
\]
for $\alpha=1$ or $2$  (where $\log^{\alpha}(x)$ stands for $(\log(x))^{\alpha}$), and that  $\widehat{p}$ estimates the ultimate proportion $p$ of uncensored observations in the tail, which turns out to be equal to
\[
p := \frac{\gamma}{\gamma_X} = \frac{\gamma_C}{\gamma_X+\gamma_C}.
\]
\zdeux

Our goal is to show that relying on usual strategies in the survival analysis literature leads to estimators of $\gamma_X$ which are often sharper than those obtained by simply dividing an estimator of $\gamma$ by the proportion  of uncensored observations. By ``usual'' strategy we mean using ``Kaplan-Meier''-like random weights : we refer to \cite{WWExtremes2014} for more detailed informations concerning the origin of the two kinds of  random weights appearing in the formulas below. As a matter of fact, we define, for any given $\alpha \geq 1$,  the following two versions of randomly weighted moments of the log relative excesses  :
\begeq{momalphaKSV}
   M^{(\alpha)}_{n,k_n}  :=  \frac 1{ n(1-\hat F_n(Z_{n-k_n,n}))} \sum_{i=1}^{k_n}
 \frac{\delta_{n-i+1,n}}{1-\hat G_n(Z^-_{n-i+1,n})} \left( \log^{\alpha} \left(  \frac{Z_{n-i+1,n}}{Z_{n-k_n,n}} \right) \, \right)
\fineq
and
\begeq{momalphaLeurg}
   \widetilde{M}^{(\alpha)}_{n,k_n}  :=  \frac 1{ n(1-\hat F_n(Z_{n-k_n,n}))} \sum_{i=1}^{k_n}
 \frac{1}{1-\hat G_n(Z^-_{n-i+1,n})} \xi_{i,n}
\fineq
where
\begeq{defxiin}
 \xi_{i,n} := i \left( \log^{\alpha} \left(  \frac{Z_{n-i+1,n}}{Z_{n-k_n,n}} \right) -  \log^{\alpha} \left(  \frac{Z_{n-i,n}}{Z_{n-k_n,n}} \right) \right)
\fineq
and $(k_n)$ is a sequence of integers satisfying, as   $n$ tends to   $+\infty$,
\begeq{kn}
 k_n \rightarrow  + \infty  \mbox{ \ and \ }  k_n= o(n).
\fineq
Above, $\hat F_n$ and $\hat G_n$ naturally denote the Kaplan-Meier estimators of $F$ and $G$, respectively, defined as follows : for $t < Z_{n,n}$,
 \[
 1-\hat F_n (t) = \prod_{Z_{i,n} \leq t} \left( \frac{n-i}{n-i+1} \right)^{\delta_{i,n}} \mbox{ and } \ \  1-\hat G_n(t) = \prod_{Z_{i,n} \leq t} \left( \frac{n-i}{n-i+1} \right)^{1-\delta_{i,n}}.
 \]

 It should be noted that these 2 weighted versions of the moments of the log-excesses defined in (1) and (2) are in fact closely related : as a matter of fact, they differ only when the maximum observation $Z_{n,n}$ is censored (when $\delta_{n,n}=1$, we have indeed $M^{(\alpha)}_{n,k_n}=\widetilde M^{(\alpha)}_{n,k_n}$ , see Proposition \ref{lemmeegaliteMnMntilde} in Section \ref{preuveTh1}). However, both versions deserve attention : firstly because in practice the last observation is often a censored one, and secondly because when they do differ, the difference is the only term involving the information contained in the maximum observation $Z_{n,n}$ (this difference is therefore non-asymptotically not negligible, although it tends to $0$ in probability, as stated in Proposition \ref{lemmeegaliteMnMntilde} in Section \ref{preuveTh1}). \zun

In section \ref{methode} below, assumptions are presented and discussed, convergence results for the weighted moments ${M}^{(\alpha)}_{n,k_n} $ and $\widetilde{M}^{(\alpha)}_{n,k_n}$ are stated, and we describe how classes of estimators of $\gamma_X$ can be deduced by combining these moments  for different values of $\alpha$. In Section \ref{simuls}, performance of these estimators will be presented on the basis of simulations. Section \ref{concl} provides some words of conclusion, Section  \ref{preuveTh1} is devoted to the proof of Theorem 1  below, and finally the Appendix includes standard (but central to our proofs) results on regularly varying functions, as well as the proofs of the different lemmas which were used in Section \ref{preuveTh1}.

\section{Results}
\label{methode}

\subsection{Assumptions}
\label{hypo}

In addition to (\ref{kn}), our results need the following minimal assumption :
\zun
\[
\mbox{(A)\hspace{0.2cm} } F\in D(H_{\gamma_X}), \hspace{0.2cm} G\in D(H_{\gamma_C})
\mbox{ \ with } \gamma_X <0\; , \gamma_C <0 \mbox{ \ and \  } x^*:=\tau_F=\tau_G.
\]
As noted earlier, this assumption implies that $H  \in D(H_{\gamma})$  with $\tau_H=x^*$ and
$$
 \gamma= \frac{\gamma_X \gamma_C}{\gamma_X+ \gamma_C} <0.
$$

\noindent If we note $U(t)=H^{\leftarrow}(1-1/t)$ the quantile function associated to $H$, then $x^*=U(\infty)$ and  $H \in D(H_{\gamma})$  is equivalent to the existence of some positive function $a$ such that
\begeq{CondLogU}
\lim_{t\rightarrow +\infty} \frac{\log U(tx)- \log U(t) }{a(t)/U(t)} = \frac{x^{\gamma}-1}{\gamma} , \ \forall x >0,
\fineq
which, since $\gamma<0$, is itself equivalent to
\begeq{CondU}
\lim_{t\rightarrow +\infty} \frac{U(\infty)- U(tx)}{U(\infty)- U(t)} = x^{\gamma}, \ \forall x >0.
\fineq
This  means that the function $U(\infty)- U$ is regularly varying (at $+ \infty$) with index $\gamma$ (see the appendix for the definition of regular variation at $+ \infty$). A reference for the equivalence of conditions (\ref{CondLogU}) and (\ref{CondU}) to (A) is \cite{deHaanFerreira2006} (respectively relation (3.5.4) and Corollary 1.2.10 there).
\zdeux

 \noindent Finally, we will need some very mild additional assumption on $(k_n)$
\zun
\[
\parbox[t]{11.8cm}{(K) \ there exists some $\delta>0$, \ or some $\delta \geq  \displaystyle{\frac{\gamma_X-\gamma_C}{\gamma_X+\gamma_C}}$    if  $\gamma_C \geq  \gamma_X$, such that
\begeq{condkn}
  -\log(k_n/n) \big/ k_n =   O(n^{-\delta})	.
\fineq
}
\]

\subsection{Asymptotic results}
\label{resultats}

Let us introduce the notation $a_{n,k} := a(n/k_n) / U(n/k_n)$ (see the previous paragraph for the definition of functions $U$ and $a$), where $a_{n,k}\rightarrow 0$ ({\it{cf}} equation (3.5.5) in \cite{deHaanFerreira2006}). In the paper, $Beta(\cdot,\cdot)$ denotes the usual Beta function, $Beta(a,b)=\int_0^1 t^{a-1}(1-t)^{b-1}\, dt$ \ $(a>0,b>0)$.
\begin{theo}
\label{consistKSV}
Under assumption (A) and conditions  $(\ref{kn})$ and (K), for any given $\alpha\geq 1$, both
$\frac{ M^{(\alpha)}_{n,k_n}}{(a_{n,k})^{\alpha}} $ and  $\frac{\widetilde{M}^{(\alpha)}_{n,k_n}}{(a_{n,k})^{\alpha}}$ converge in probability, as $n$ tends to $\infty$, to
\[
|\gamma_X|^{-1} |\gamma|^{-\alpha} Beta(|\gamma_X|^{-1} \ ; \  \alpha+1).
\]
\end{theo}

The following corollary states the consistency of our two  different adaptations of the Moment estimator to this censored framework.

\begin{cor}
\label{Mom}
Under conditions of Theorem \ref{consistKSV}, as $n\rightarrow\infty$,
\[
\widehat{\gamma}_{n,Mom} :=  M^{(1)}_{n,k_n}  + 1 -   \frac 12 \left(1- \frac{(M^{(1)}_{n,k_n})^2}{M^{(2)}_{n,k_n}} \right)^{-1}  \ \limprob \ \gamma_X
\]
and
\[
\widetilde{\gamma}_{n,Mom} :=  \widetilde{M}^{(1)}_{n,k_n}  + 1 -   \frac 12 \left(1- \frac{(\widetilde{M}^{(1)}_{n,k_n})^2}{\widetilde{M}^{(2)}_{n,k_n}} \right)^{-1}  \ \limprob \ \gamma_X.
\]
\end{cor}

In fact, by using the elementary properties of the Beta function, the weighted moments $M_n^{(\alpha)}$ or $\widetilde{M}_n^{(\alpha)}$ can be combined in different ways, leading to the definition of two different classes of consistent estimators of  $\gamma_X$, parametrized by $\alpha\geq 1$ (proofs of the 3 corollaries are easy and omitted). In the next section, we study their finite sample performance.
\begin{cor}
\label{alphaMom1}
Under  conditions of Theorem \ref{consistKSV}, as $n\ra\infty$,
\[
\widehat{\gamma}^{(\alpha)}_{n,1} := \left( V^{-1}_{n,\alpha} + \alpha +1 \right)^{-1} \ \limprob \ \gamma_X
\]
and
\[
\widetilde{\gamma}^{(\alpha)}_{n,1} := \left( \widetilde{V}^{-1}_{n,\alpha} + \alpha +1  \right)^{-1} \ \limprob \ \gamma_X
\]
where
\[
V_{n,\alpha} := 1- \frac{\alpha + 2}{\alpha + 1} \ \frac{(M_n^{(\alpha +1)})^2}{M_n^{(\alpha)} M_n^{(\alpha +2)}}
\makebox[1.5cm][c]{and}
\widetilde{V}_{n,\alpha} := 1- \frac{\alpha + 2}{\alpha + 1} \ \frac{(\widetilde{M}_n^{(\alpha +1)})^2}{\widetilde{M}_n^{(\alpha)} \widetilde{M}_n^{(\alpha +2)}}.
\ztrois
\]
\end{cor}

\begin{cor}
\label{alphaMom2}
Under  conditions of Theorem \ref{consistKSV}, as $n\ra\infty$,
\[
\widehat{\gamma}^{(\alpha)}_{n,2} := \frac{1-(\alpha+1) R_{n,\alpha}}{(\alpha+1) (1-R_{n,\alpha})} \ \limprob \ \gamma_X
\]
and
\[
\widetilde{\gamma}^{(\alpha)}_{n,2} := \frac{1-(\alpha+1) \widetilde{R}_{n,\alpha}}{(\alpha+1) (1-\widetilde{R}_{n,\alpha})} \ \limprob \ \gamma_X ,
\]
where
\[
R_{n,\alpha} :=  \frac{M_n^{(1)} M_n^{(\alpha)}}{M_n^{(\alpha+1)}}
\makebox[1.5cm][c]{and}
\widetilde{R}_{n,\alpha} := \frac{\widetilde{M}_n^{(1)} \widetilde{M}_n^{(\alpha)}}{\widetilde{M}_n^{(\alpha+1)}} .
\ztrois
\]
\end{cor}

\begin{rmk}
It is straightforward to see that $\hat\gamma^{(\alpha)}_{n,2}$ with $\alpha=1$ equals $1-\frac 1 2 (1-R_{n,1})^{-1}$, which is very close to $\hat{\gamma}_{n,Mom}$, since $M^{(1)}_{n,k_n}\rightarrow 0$ in our finite endpoint framework.
\end{rmk}

\begin{rmk}
\label{Momclassique}
If ${\cal M}^{(\alpha)}_{n,k_n}$ denotes the unweighted moments defined in the introduction, it can be proved that under $(A)$ and $(\ref{kn})$, for $\alpha\geq 1$,
\[
\frac{ {\cal M}^{(\alpha)}_{n,k_n}}{(a_{n,k})^{\alpha}}  \ \limprob \ |\gamma|^{-\alpha-1} Beta(|\gamma|^{-1} \ ; \  \alpha+1) \, .
\]
Therefore, it is easy to check that combining those moments as described in Corollaries \ref{alphaMom1} and \ref{alphaMom2} leads to consistent estimators of $\gamma$, and thus dividing the latter  by $\hat p$ (defined in the introduction) leads to 2 classes of consistent estimators $\widecheck\gamma^{(\alpha)}_{n,1}$ and $\widecheck\gamma^{(\alpha)}_{n,2}$ of $\gamma_X$. We also define $\widecheck\gamma_{n,Mom}$ as the estimator of $\gamma_X$ obtained by dividing the classical Moment estimator of $\gamma$ by the proportion $\hat p$. A finite-sample comparison of those  estimators with our new competitors is presented  in the following section.
\end{rmk}

\begin{rmk}
\label{casGammaPositif}
Note that the combination of moments proposed in Corollaries 2 and 3 become inadequate in the framework of a positive extreme value index: it can be indeed proved that, in this framework, the combinations $\widehat\gamma^{(\alpha)}_{n,j}$ and $\widetilde\gamma^{(\alpha)}_{n,j}$, for $j=1$ or $2$, converge in probability to zero, by proving that $M^{(\alpha)}_{n,k_n}$ and $\widetilde M^{(\alpha)}_{n,k_n}$ converge in probability to $\gamma_X^{\alpha} \Gamma(\alpha+1)$ (in the complete data case, this result is known for ${\cal M}^{(\alpha)}_{n,k_n}$, see \cite{Segers2001}). This could suggest, in the positive index case, the definition of estimators of $\gamma_X$ which would be equal to $\hat\gamma^{(\alpha)}_{n,j}$ or $\tilde\gamma^{(\alpha)}_{n,j}$ (for $j=1$ or $2$) plus a "censored version" of the Hill estimator (in the same spirit as the definition of the Moment estimator, which equals the Hill estimator plus a term converging to $0$ in the positive index case).
\end{rmk}

\section{Finite sample behavior}
\label{simuls}

The goal of this Section  is to present our results concerning the finite sample performances of our new estimators of the extreme value index  in presence of random censorship, presented in Corollaries \ref{Mom},  \ref{alphaMom1} and  \ref{alphaMom2}. In each case considered,  $2000$ random samples of size $n=500$ were generated, and the median bias and mean squared error (MSE) of the different estimators of $\gamma$ were plotted against the number $k_n$ of excesses used.
\zdeux

A great variety of situations can be (and has been) considered in our simulation study : various values of $\gamma_X$ and $\gamma_C$ (and therefore various censoring rates in the tail), various families of underlying distributions (Reverse Burr, generalized Pareto, Beta), and choice of the value of $\alpha$. It is impossible to illustrate here the different possible combinations of these features : we will therefore try to draw some general conclusions from the many different situations we have observed, and provide a partial illustration with 3 particular cases.

Concerning the choice  of the tuning parameter $\alpha$, we did not find a value which seemed preferable in every situation : nonetheless, in general, for small values of $k_n$, a value of $\alpha$ around 1 or 2 yields better MSE, whereas for high values of $k_n$, the MSE is lower for values of $\alpha$ greater than 2. We decided not to include this preliminary  study in this article, and chose (almost arbitrarily) the value $\alpha=2$ in all our subsequent simulations.
\zdeux

Let us now settle the vocabulary used in this section. We will call {\it Moment estimators} the estimators $\widehat\gamma_{n,Mom}$ and $\widetilde\gamma_{n,Mom}$ appearing in Corollary \ref{Mom}, as well as the estimator $\widecheck\gamma_{n,Mom}$ introduced in Remark \ref{Momclassique} above. We will call {\it type 1 (resp. type 2) estimators} the estimators $\widehat\gamma^{(\alpha)}_{n,1}$ and $\widetilde\gamma^{(\alpha)}_{n,1}$ (resp. $\widehat\gamma^{(\alpha)}_{n,2}$ and $\widetilde\gamma^{(\alpha)}_{n,2}$) appearing in Corollary \ref{alphaMom1} (resp. \ref{alphaMom2}), as well as the estimator $\widecheck\gamma^{(\alpha)}_{n,1}$ (resp. $\widecheck\gamma^{(\alpha)}_{n,2}$) introduced in Remark \ref{Momclassique}.
\zun

\begin{figure}[htp]
\centering
\includegraphics[height=7.0cm,width=.9\textwidth]{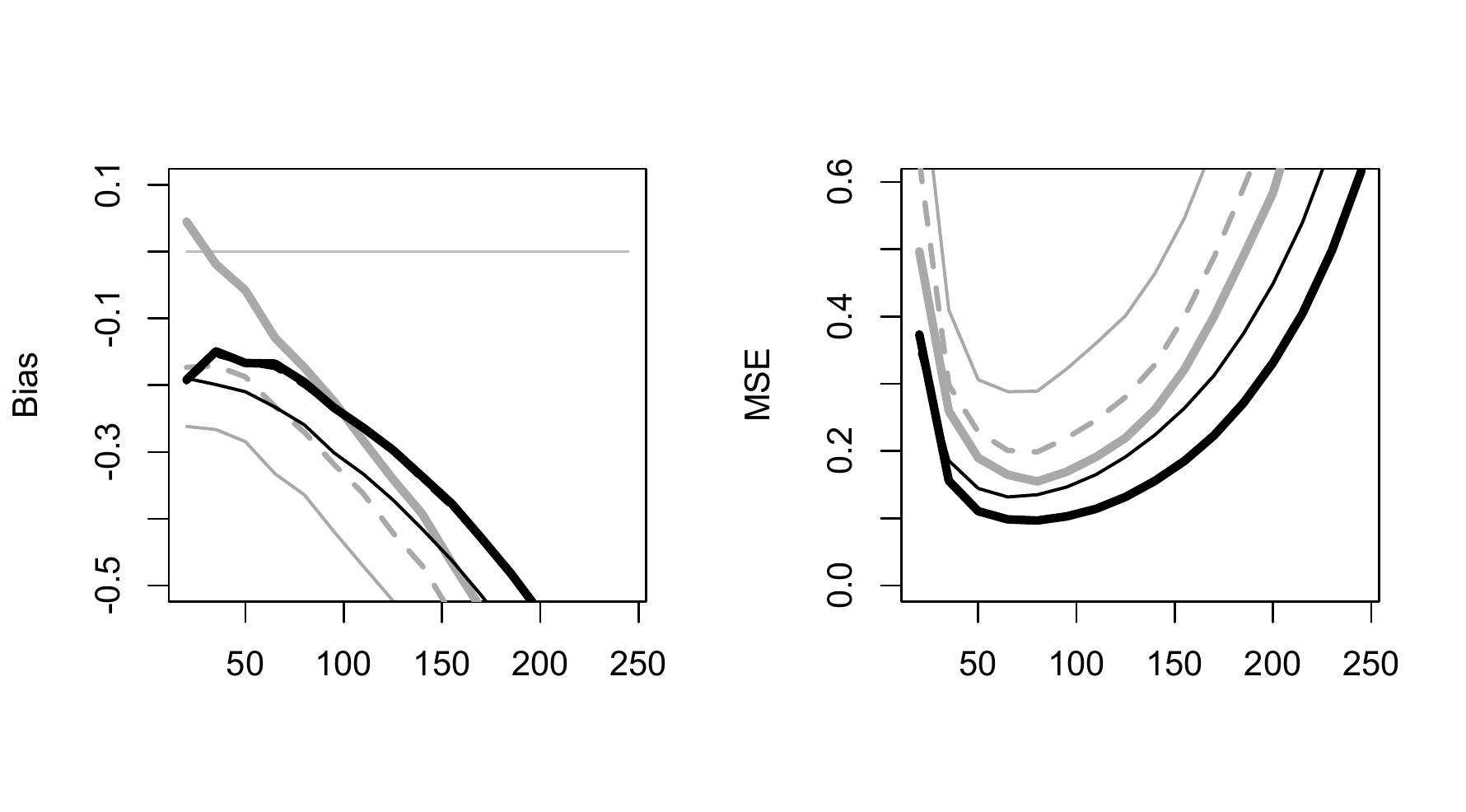}
\caption{Comparison between $\widehat{\gamma}^{(2)}_{n,1}$ (thick black), $\widetilde{\gamma}^{(2)}_{n,1}$ (dashed black), $\widecheck{\gamma}^{(2)}_{n,1}$ (thin black), $\widehat{\gamma}_{n,Mom}$ (thick grey), $\widetilde{\gamma}_{n,Mom}$ (dashed grey) and $\widecheck{\gamma}_{n,Mom}$ (thin grey) for a RevBurr$(1,1,1,10)$  censored by a RevBurr$(10,2/3,1,10)$ ($\gamma_X=-1 >\gamma_C=-3/2$, p=2/5, weak censoring)}
\end{figure}

We will also consider names for the different methods : the KM method (for Kaplan-Meier-like weights, appearing in the definition of $M^{(\alpha)}_{n,k_n}$), leading to $\widehat\gamma$ estimators, the L method (for Leurgans-like weights) leading to $\widetilde\gamma$ estimators (the name comes from the mathematician Sue Leurgans who inspired the weights, see \cite{WWExtremes2014} for details and a reference), and the EFG method (for constant weighting by $\hat p$), leading to $\widecheck\gamma$ estimators (the names comes from the initials of the authors of \cite{Einmahl2008}).
\zdeux

There are two main questions addressed in this empirical study : is one of the 3 methods preferable to the others (and in which conditions) and is there a better choice for the type of estimator (type 1 , type 2, or classical Moment estimator) ? Unsurprisingly, after our intensive simulation study, we may say that the answer is no for the 2 questions, if an overall superiority is looked for. However, we can make some partial comments concerning the choice of the method and of the estimator type, whether the censoring is strong or weak, or the value of $|\gamma|$ is small or not.
\zdeux

Note first that, if the censoring rate $1-p$ in the tail is very low (say  lower than $10\%$), we observed that there was not much difference between the 3 methods (KM, L, EFG), and that  it was just a question of choosing between type 1, type 2, and moment estimator. This is why, in the following, we only consider cases where the censoring rate $1-p=\frac{\gamma_X}{\gamma_X+\gamma_C}$ is larger than 1/4, and talk about strong censoring in the tail when this rate is greater than $1/2$ ({\it i.e.} $\gamma_X\leq\gamma_C$),  and weak censoring otherwise (when $\gamma_X>\gamma_C$).
\zdeux

For ``high'' values of $\gamma_X$, {\it i.e.} lower than $-1/2$, we have most of the time observed better performance of the KM and L methods with respect to the EFG method, in strong or weak censoring frameworks.  In this context, the type 1 estimators are generally preferable to the type 2 estimators, and comparable or preferable to the moment estimator.
 \zdeux

\begin{figure}[htp]
\centering
\includegraphics[height=7.0cm,width=.9\textwidth]{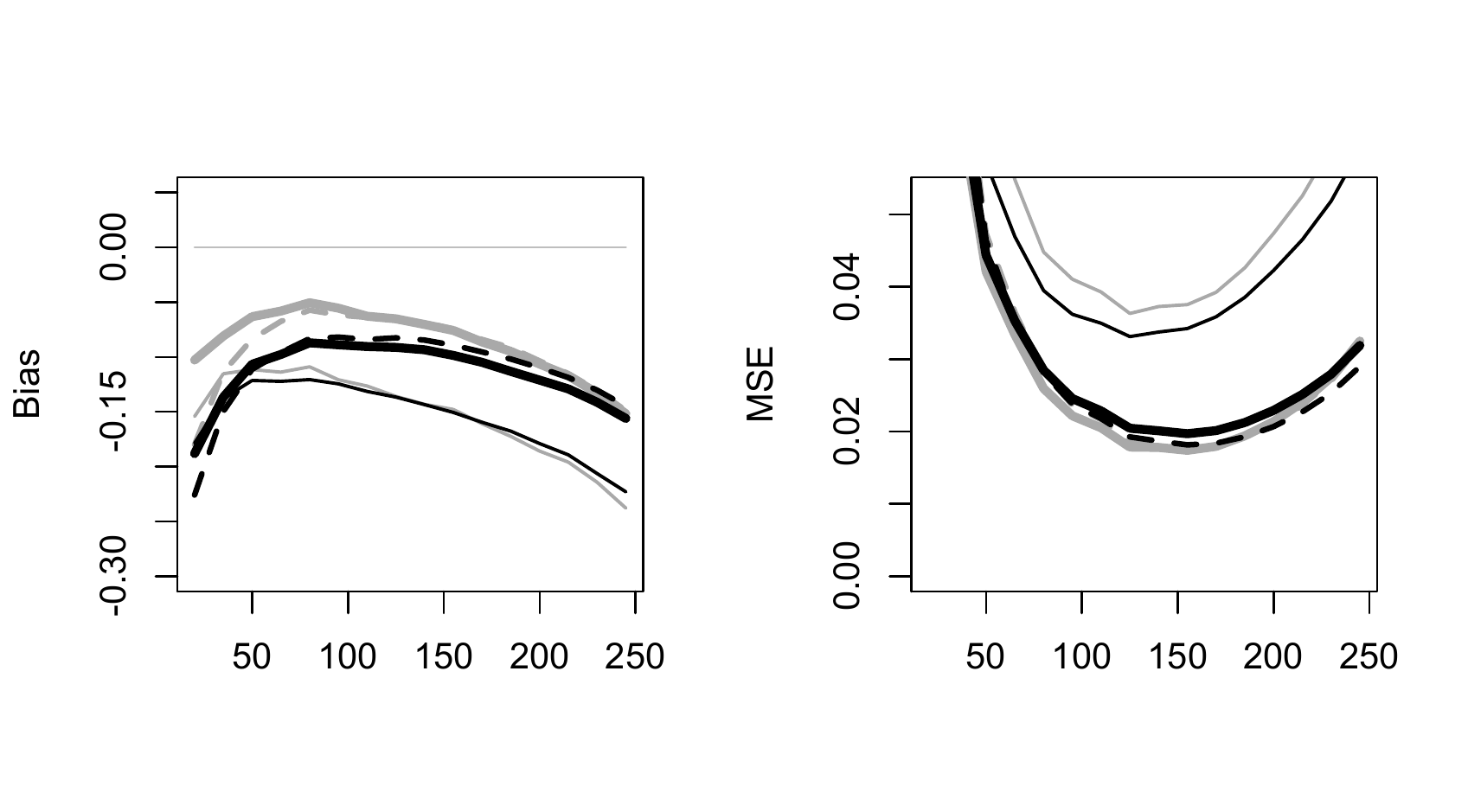}
\caption{Comparison between $\widehat{\gamma}^{(2)}_{n,2}$ (thick black), $\widetilde{\gamma}^{(2)}_{n,2}$ (dashed black), $\widecheck{\gamma}^{(2)}_{n,2}$ (thin black), $\widehat{\gamma}_{n,Mom}$ (thick grey), $\widetilde{\gamma}_{n,Mom}$ (dashed grey) and $\widecheck{\gamma}_{n,Mom}$ (thin grey)   for a RevBurr$(1,8,1/2,10)$  censored by a RevBurr$(10,4,1/2,10)$ ($\gamma_X=-1/4 >\gamma_C=-1/2$, p=1/3, weak censoring)}
\end{figure}

For values of $\gamma_X$ between $-1/2$ and $0$ (sometimes called the ``regular'' case, and which is the most frequently encountered in practice), there exists a great variety of situations. We observed that the moment estimators were generally better than the type 2 estimators, which were themselves generally better than the type 1 ones. Concerning the choice of the method, for the moment estimator, it seems difficult to suggest a particular one, between the KM, L, and EFG methods (even though in many cases, at least one among the KM and L methods was better than the EFG method).  Concerning the inferiority of types 1 and 2 versus the moment estimator, it should be noted that it is mainly due to the bias, which contributes the most to the MSE (in fact, we clearly noticed that the variances of the types 1 and 2, for $\alpha=2$, are almost always lower than the variance of the moment estimator).
\ztrois

The 3 particular situations we chose as illustrations of the comments above  involve the Reverse Burr class of distributions  $RevBurr(\beta,\tau,\lambda,x^*)$  (with $\beta,\tau,\lambda>0$) : its survival function is
\[
\bP(X>x) = (1+\beta^{-1} (x^*-x)^{-\tau})^{-\lambda},
\]
 and its extreme value index  is $-1/(\lambda \tau)$.
\zdeux

In Figure 1, the value of $\gamma_X$ is lower than $-1/2$, and therefore, as motivated above, for readability purposes we only kept the type 1 estimators on the graph, whereas for the other two figures, the value of $\gamma_X$ is between $-1/2$ and $0$ and we therefore only kept the type 2 estimator illustrated. Remind here that these 3 examples are only 3 particular cases of the numerous combinations of features we have considered in our simulation study.

\begin{figure}[htp]
\centering
\includegraphics[height=7.0cm,width=.9\textwidth]{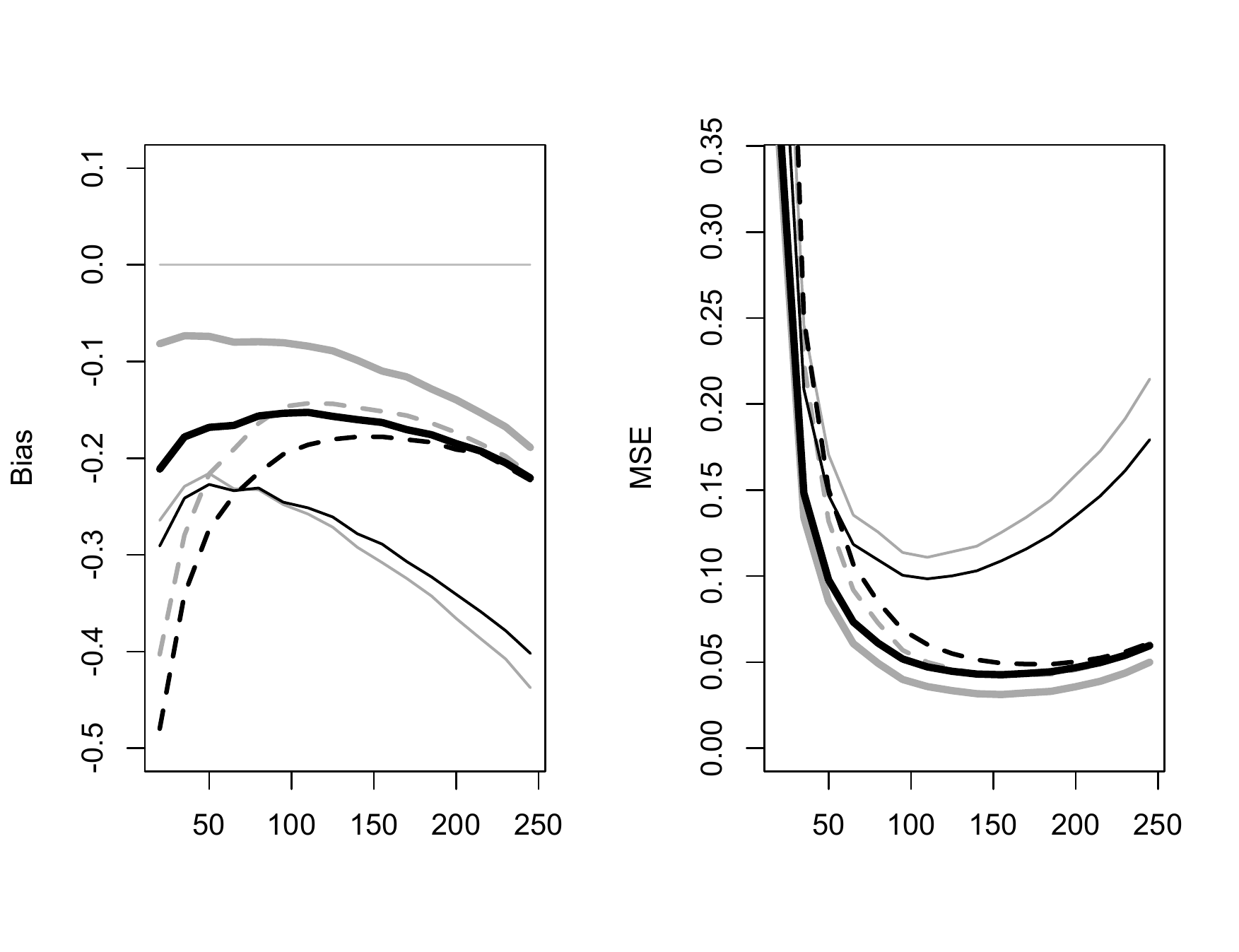}
\caption{Comparison between $\widehat{\gamma}^{(2)}_{n,2}$ (thick black), $\widetilde{\gamma}^{(2)}_{n,2}$ (dashed black), $\widecheck{\gamma}^{(2)}_{n,2}$ (thin black), $\widehat{\gamma}_{n,Mom}$ (thick grey), $\widetilde{\gamma}_{n,Mom}$ (dashed grey) and $\widecheck{\gamma}_{n,Mom}$ (thin grey)   for a RevBurr$(10,8,1/2,10)$  censored by a RevBurr$(10,5,1,10)$ ($\gamma_X=-1/4<\gamma_C=-1/5$, p=5/9, strong censoring)}
\end{figure}


\section{Conclusion}
\label{concl}

In this paper, we applied the methodology introduced in \cite{WWExtremes2014}  to define weighted versions of the moments of relative excesses, and consequently construct new estimators of the extreme value index  for randomly-censored data with distributions in the Weibull domain of attraction. We proposed, in particular, a new adaptation of the famous Moment estimator.   Our intensive simulation study shows that  the proposed estimators are  competitive even if, in many cases,  the bias would need to be reduced.  A future possible work would be to exploit our weighting methodology in order to estimate other parameters of the tail (for reducing  the bias, for example) as well as extreme quantiles.  The asymptotic normality remains a question to be addressed (difficulties come from the control of the Kaplan-Meier estimates in the tail).

\newpage

\section{Proof of Theorem \ref{consistKSV}} \label{preuveTh1} \zdeux

Before proceeding to the proof of Theorem \ref{consistKSV}, we state the following Proposition which explains the link between  our two proposals of weighted moments.
\begin{prop}
\label{lemmeegaliteMnMntilde}
\begitem
\item[$(i)$] For any $\alpha \geq 1$,
 $\widetilde{M}^{(\alpha)}_{n,k_n} =  M^{(\alpha)}_{n,k_n} + (1-\delta_{n,n}) D^{(\alpha)}_{n,k_n}$,
where
\[
 D^{(\alpha)}_{n,k_n} = \frac 1{ n(1-\hat F_n(Z_{n-k_n,n})) (1-\hat G_n(Z^-_{n,n}))} \  \log^{\alpha} \left(  \frac{Z_{n,n}}{Z_{n-k_n,n}} \right).
\]
\item[$(ii)$] Under the same assumptions as Theorem \ref{consistKSV}, for any $\alpha \geq 1$, we have $D^{(\alpha)}_{n,k_n} = o_{\bP}(a^{\alpha}_{n,k})$.
\finit
\end{prop}
\noindent According to this proposition,  the validity  of  Theorem \ref{consistKSV} for $M^{(\alpha)}_{n,k_n}$ is a consequence of  Theorem \ref{consistKSV} for $ \widetilde{M}^{(\alpha)}_{n,k_n}$.  Proposition \ref{lemmeegaliteMnMntilde}  will be proved at the end of  the appendix. \zdeux

 \noindent  We now  need to state the following technical  Lemmas, which will be proved in the appendix.  \mbox{The notation $\xi_{i,n}$ was introduced in (\ref{defxiin})} and we note
 \[
 \tilde{Z}_{i,n} := \frac{x^* - Z_{n-i+1,n}}{x^* - Z_{n-k_n,n}} .
 \]

 \begin{lem}
\label{lemmeXin}
Let $\alpha \geq 1$ and $u_i= u_{i,n} = \frac{i}{k_n+1}$. Under assumptions (A) and (\ref{kn}):
\zun\\
$(i)\;$ if $0 < a < 1$ then
\[
\frac{1}{k_n} \sum_{i=1}^{k_n} u_i^{-a} \frac{\xi_{i,n}}{a_{n,k}^{\alpha}}  \ \limprob \ (1-a) |\gamma|^{-\alpha-1} Beta\left(\frac{1-a}{|\gamma|} \ ; \  \alpha +1 \right) \, ,
\]
$(ii)\;$ if $a >1$,  then for any $\delta'>0$,
 \[
 \frac{1}{k_n^{a+\delta'}} \sum_{i=1}^{k_n}  u_i^{-a} \frac{\xi_{i,n}}{a_{n,k}^{\alpha}}  \ \limprob \  0.
 \]
\end{lem}

\begin{lem}
\label{lemmeBin}
For any given positive exponents $\theta$ and $\theta'>0$, there exist constants $c>1$, $c'<1$ both arbitrarily close to 1, and
$a_+ >0$, $a_->0$, arbitrarily close to $\gamma\theta$ and to $\gamma\theta'$ respectively,
such that
 \[
 \lim_{n \rightarrow \infty} \bP \left( \max_{i \leq k_n} \frac{ \tilde{Z}_{i,n}^{\theta}  }{c u_i^{-a_+}} > 1\right)
 =  \lim_{n \rightarrow \infty} \bP \left( \min_{i \leq k_n} \frac{ \tilde{Z}_{i,n}^{\theta'} }{c' u_i^{-a_-}} < 1 \right)
 = 0.
 \]
\end{lem}
\zdeux

\noindent We now proceed  to the proof of Theorem \ref{consistKSV} (for $ \widetilde{M}^{(\alpha)}_{n,k_n}$), which has structural similarities with the proof of Theorem 2 in \cite{WWExtremes2014}. We shall refer to the latter when necessary. We  have the decomposition
\[
 \widetilde{M}^{(\alpha)}_{n,k_n} = A_n  ( \overline{W}_n +  R_n) \; ,
\]
where
\[
 \begar{lcl}
  A_n := \frac{1-\widehat G_n(Z_{n-k_n,n})}{1-G(Z_{n-k_n,n})}  & , & \hspace{0.3cm} \overline{W}_n := \frac{1}{k_n} \sum_{i=1}^{k_n}    W_{i,n},    \zdeux \\
  R_n := \frac{1}{k_n} \sum_{i=1}^{k_n}  (C_{i,n}-1)   W_{i,n}  \hspace{0.3cm} {} & , &  \hspace{0.3cm} C_{i,n}  : = \frac{1- G(Z_{n-i+1,n})}{1-\widehat G_n(Z^-_{n-i+1,n})},
 \finar
\]
and
\[
   W_{i,n}   :=  \frac{1-G(Z_{n-k_n,n})}{1- G(Z_{n-i+1,n})} \xi_{i,n}.
\]
 Since $A_n \limprob 1$ as $n \rightarrow\infty$ (see Theorem 2 in \cite{Csorgo96}), we need to prove that $R_n = o_{\bP}(a_{n,k}^{\alpha})$ and
 \begeq{convWbarn}
 \frac{\overline{W}_n}{a_{n,k}^{\alpha}} \limprob l_{\alpha} ,
 \fineq
where $l_{\alpha}$ denote the limit in the statement of Theorem \ref{consistKSV}.

\subsection{Proof of  $\overline{W}_n / a_{n,k}^{\alpha} \limprob l_{\alpha}$}
\zun
Since $G\in D(H_{\gamma_C})$, with $\gamma_C <0$, is equivalent to $t \rightarrow 1-G(x^*-t)$ being regularly varying at $0$ with index $-1/\gamma_C$  (see the appendix for the definition of regular variation at $0$), the bounds $(\ref{Potter1surf})$ in Corollary \ref{BornesPotter0} (in the appendix) applied to $f=G$, $x=(x^*-Z_{n-i+1,n})/(x^*-Z_{n-k_n,n})$ and $t=x^*-Z_{n-k_n,n}$, yield,  for $\epsilon >0$, $n$ sufficiently large  and every $1 \leq i \leq k_n$,
 \begeq{bornesWin}
 (1-\epsilon)  \ \xi_{i,n} \  \tilde{Z}_{i,n}^{\gamma_C^{-1} + \epsilon}  \leq W_{i,n} \leq (1+\epsilon)   \ \xi_{i,n} \ \tilde{Z}_{i,n}^{\gamma_C^{-1} - \epsilon}.
 \fineq
Let $\eta >0$.  We first write, for $\epsilon$ sufficiently small,
\begin{eqnarray*}
\bP(\; \overline{W}_n/a_{n,k}^{\alpha}  - \l_{\alpha} > \eta \; ) \leq \bP( \; k_n^{-1}  \textstyle\sum_{i=1}^{k_n} \tilde{Z}_{i,n} ^{\gamma_C^{-1}-\epsilon} \frac{\xi_{in}}{a_{n,k}^{\alpha}} -  \l_{\alpha} > \frac{\eta}{2} \; ), \\
\bP(\; \l_{\alpha} - \overline{W}_n/a_{n,k}^{\alpha}  > \eta \; ) \leq \bP(\; l_{\alpha} - k_n^{-1}  \textstyle \sum_{i=1}^{k_n}  \tilde{Z}_{i,n} ^{\gamma_C^{-1}+\epsilon} \frac{\xi_{in}}{a_{n,k}^{\alpha}}> \frac{\eta}{2} \; ).
\end{eqnarray*}
Let us now consider constants $c >1$ and $c'<1$, both arbitrary close to $1$, and $a_+ >0$ and $a_-> 0$ both arbitrary close to $\gamma/\gamma_C$ . These constants come from the application of Lemma \ref{lemmeBin} above with $\theta=\gamma_C^{-1}-\epsilon$ and $\theta'=\gamma_C^{-1}+\epsilon$.  Using positivity of $ \xi_{in}$, it comes
\begin{eqnarray*}
\nonumber
\bP\left(\;\frac{\overline{W}_n}{a_{n,k}^{\alpha}} - \l_{\alpha} > \eta\;\right)
& \leq & \bP\left(\max_{i \leq k_n} \frac{ \tilde{Z}_{i,n} ^{\gamma_C^{-1}-\epsilon}  }{c u_i^{-a_+}} > 1\right) \\
& & \hspace*{0.6cm} + \;
  \bP \left( c \ k_n^{-1} \sum_{i=1}^{k_n} u_i^{-a_+} \frac{\xi_{in}}{a_{n,k}^{\alpha}} - \l_{\alpha} > \frac{\eta}{2} \right)
\\
\nonumber
\bP\left(\;\l_{\alpha} -\frac{\overline{W}_n}{a_{n,k}^{\alpha}}  > \eta\;\right)
& \leq & \bP\left(\min_{i \leq k_n} \frac{  \tilde{Z}_{i,n} ^{\gamma_C^{-1}+\epsilon}  }{c' u_i^{-a_-}} < 1 \right)  \\
& & \hspace*{0.6cm} + \;
  \bP \left( \l_{\alpha} - c' \ k_n^{-1} \sum_{i=1}^{k_n} u_i^{-a_-} \frac{\xi_{in}}{a_{n,k}^{\alpha}} > \frac{\eta}{2} \right)
\end{eqnarray*}
where $u_i=u_{i,n} = i/(k_n+1)$ for $1 \leq i \leq k_n$. If we call $l_{\alpha,a}$ the limit in the statement of Lemma \ref{lemmeXin} $(i)$, and if we apply Lemma \ref{lemmeBin} as indicated previously, we have
\[
 \begar{rcl}
  {\displaystyle \limsup_{n\ra\infty}} \; \bP\left(\; \frac{\overline{W}_n}{a_{n,k}^{\alpha}} - l_{\alpha} > \eta  \;\right) & \leq &
   {\displaystyle \limsup_{n\ra\infty}} \;  \bP \left(  k_n^{-1} \sum_{i=1}^{k_n} u_i^{-a_+} \frac{\xi_{in}}{a_{n,k}^{\alpha}} - l_{\alpha,a_+} > \frac 1 c (\frac{\eta}{2}+l_{\alpha}) -  l_{\alpha,a_+} \right)
  \zun\\
 {\displaystyle \limsup_{n\ra\infty}} \;\bP\left(\; l_{\alpha} -\frac{\overline{W}_n}{a_{n,k}^{\alpha}}  > \eta\;\right) & \leq &
 {\displaystyle \limsup_{n\ra\infty}} \;\bP \left(  l_{\alpha,a_-} \! - k_n^{-1} \sum_{i=1}^{k_n} u_i^{-a_-} \frac{\xi_{in}}{a_{n,k}^{\alpha}} > \frac 1 {c'} (\frac{\eta}{2}-l_{\alpha}) +  l_{\alpha,a_-} \! \right)
 \finar
\]
\zun
Since $l_{\alpha}=l_{\alpha,\gamma/\gamma_C}$, and both $a_+$ and $a_-$ are arbitrary close to $\gamma/\gamma_C<1$, it is easy to see that $(\ref{convWbarn})$ comes from the application of Lemma \ref{lemmeXin} $(i)$ to $a=a_+$ and $a=a_-$.
\ztrois

\subsection{Proof of $R_n=o_{\bP}(a_{n,k}^{\alpha})$}
\label{PreuveRnKSV}
Let us use the same decomposition as in the proof of the negligibility of the term $R_n$ in  \cite{WWExtremes2014} (see subsection 5.1.2 there). In other words, we define, for some $\delta'>0$,
 \[
   \tilde{C} (t) : = \int^t_0 \ \frac{dG(x)}{(1-G(x))^2 (1-F(x))} \makebox[1.5cm][c]{ and }  h_{i,n} :=(\tilde{C} (Z_{n-i+1,n}) )^{-\frac12 -\delta'},
 \]
and we readily have  $| R_n | \leq   T^1_n  T^2_n$,  where
 \[
  T^1_n :=   \sup_{1\leq i \leq k_n}  \sqrt{n} |h_{i,n} (C_{i,n}-1)|
  \makebox[1.5cm][c]{ and }
  T^2_n :=\frac{1}{k_n} \sum_{i=1}^{k_n}  W_{i,n} h^{-1}_{i,n} n^{-\frac12}.
 \]
Using sharp results of the survival analysis literature, we have already proved in \cite{WWExtremes2014} that $ T^1_n= O_{\bP}(1)$. It remains to prove that
$$
T^2_n / a_{n,k}^{\alpha}= o_{\bP}(1).
$$
First, from the definition of $h_{in}$ and $\tilde C$, since $(1-H)=(1-F)(1-G)$ we clearly have
\[
 h^{-1}_{in} < \left( \frac{-\log(1-G(Z_{n-i+1,n}))}{1-H(Z_{n-i+1,n})} \right)^{{\frac12+\delta'}}.
\]
Moreover, under assumption (A), $1-H(x^*- \cdot)$ is regularly varying at zero with index $-1/\gamma$ and $-\log(1-G(x^*- \cdot))$ is slowly varying at $0$ :  therefore,  the application of $(\ref{bornesWin})$, as well as  bound $(\ref{Potterf})$ to $-\log(1-G(x^*- \cdot))$  and bound (\ref{Potter1surf}) to $f=G$ and $f=H$, implies that , for $n$ sufficiently larger,  $T^2_n \leq 4 P_n Q_n$  where
\begin{eqnarray*}
   P_n & := &  n^{-\frac12} \left( \frac{-\log(1-G(Z_{n-k_n,n}))}{1-H(Z_{n-k_n,n})} \right)^{{\frac12+\delta'}} \\
   Q_n & := &  \frac{1}{k_n} \sum_{i=1}^{k_n}   \xi_{in} \  \tilde Z_{i,n}^{\beta},
 \end{eqnarray*}
where $\beta=(2\gamma)^{-1} +  \gamma_C^{-1} - \epsilon'$, for some $\epsilon'>0$.  \zun \\
We thus need to prove that $P_nQ_n / a_{n,k}^{\alpha}=o_{\bP}(1)$. \zun \\

Let $\eta >0$ and consider constants $c >1$ arbitrarily close to $1$ and $a_+ >0$  arbitrarily close to $\gamma\beta=\frac 1 2 + \frac{\gamma}{\gamma_C}-\gamma\epsilon'  > \frac 1 2 + \frac{\gamma}{\gamma_C}$. We have,
\begeq{borneT2nLeurg}
\bP \left( \frac{P_nQ_n}{a_{n,k}^{\alpha}}>\eta \right) \leq  \bP\left(\max_{i \leq k_n} \frac{ \tilde{Z}_{i,n} ^{\beta}  }{c u_i^{-a_+}} > 1\right) +
  \bP \left( P_n \ k_n^{-1} \sum_{i=1}^{k_n} u_i^{-a_+} \frac{\xi_{in}}{a_{n,k}^{\alpha} } > \frac{\eta}{c} \right).
\fineq

First, Lemma \ref{lemmeBin}  is applied with $\theta= \beta$ and thus the first term  of the right-hand side of  $(\ref{borneT2nLeurg})$ tends to $0$.  Next, Lemma \ref{lemmeXin} is applied with $a=a_+$ : we thus need to distinguish  the case $\gamma_C < \gamma_X$ (for which $\gamma\beta<1$) from the case $\gamma_C \geq \gamma_X$ (for which $\gamma\beta>1$ when $\epsilon'>0$ gets small).
\begitem
 \item[$(i)$] Case  $\gamma_C < \gamma_X$
\zun\\
First of all,  assumption (K) implies that  $P_n = o_{\bP}(1)$ (see relation $(20)$ in \cite{WWExtremes2014}). Since $a_+ <1$ in this case, Lemma \ref{lemmeXin} $(ii)$ implies that $\frac{1}{k_n}  \sum_{i=1}^{k_n} u_i^{-a_+} \frac{\xi_{in}}{a_{n,k}^{\alpha}}  = O_{\bP}(1)$ and consequently the second term  of the right hand-side of ($\ref{borneT2nLeurg}$) tends to $0$.
\zdeux

\item[$(ii)$] Case  $\gamma_C \geq  \gamma_X$
\zun\\
In this case $a_+>1$,  therefore  \mbox{Lemma  \ref{lemmeXin} $(ii)$ implies that, for any given $\delta' >0$,}
 $k_n^{-(a_+ + \delta')}  \sum_{i=1}^{k_n} u_i^{-a_+} \frac{\xi_{in}}{a_{n,k}^{\alpha}}  = o_{\bP}(1)$. Moreover, assumption (K) implies that, for $\delta'>0$ small enough and $a_+$ sufficiently close to $1/2+\gamma/\gamma_C$,  we have $k_n^{a_+ + \delta' -1}P_n=O_{\bP}(1)$. Hence, the second term  of the right hand-side of ($\ref{borneT2nLeurg}$) also tends to $0$.
\hfill$\diamond$
\finit

 \section{Appendix}
\label{annexe}

 \subsection{Regular variation and Potter-type bounds}

\begin{defi} \label{defRV}
An ultimately positive function $f$ : $\bR^+  \rightarrow \bR$  is {\it regularly varying} (at infinity) with index  $\alpha \in \bR$, if
\[
\lim_{t \rightarrow + \infty} \frac{f(tx)}{f(t)} = x^{\alpha} \hspace{0.2cm}  (\forall x >0).
\]
This is noted $f \in RV_{\alpha}$.  If $\alpha=0$, $f$ is said to be slowly varying.
\end{defi}
\begin{rmk}
Regular  variation (and slow variation) can be defined at zero as well. A function $f$ is said to be regularly varying at zero with index $\alpha$ if the function $x \rightarrow  f(1/x)$  is regularly varying at infinity, with index $-\alpha$.
\end{rmk}

\begin{prop} \label{PotterBounds} (See  \cite{deHaanFerreira2006}  Proposition B.1.9) \\
Suppose $f \in RV_{\alpha}$. If $x>0$ and $\delta_1,\delta_2  >0$ are given, then there exists $t_0=t_0(\delta_1,\delta_2)$ such that for any $t\geq t_0$ satisfying $tx \geq t_0$, we have
\[
(1-\delta_1) x^{\alpha} \min(x^{\delta_2}, x^{-\delta_2})  < \frac{f(tx)}{f(t)} < (1+\delta_1) x^{\alpha} \max(x^{\delta_2}, x^{-\delta_2}) .
\]
If $x < 1$ and $\epsilon >0$, then there exists $t_0=t_0(\epsilon)$ such that for every $t\geq t_0$,
\[
(1-\epsilon) x^{\alpha+\epsilon} < \frac{f(tx)}{f(t)} < (1+\epsilon) x^{\alpha-\epsilon}
\]
and if $x \geq  1$ ,
\[
(1-\epsilon) x^{\alpha-\epsilon} < \frac{f(tx)}{f(t)} < (1+\epsilon) x^{\alpha+\epsilon} .
\]
\end{prop}

\begin{cor} \label{BornesPotter0}
If $f$ is a positive function with end-point $x^*$, such that $t \rightarrow 1-f(x^*-t)$ is regularly varying at $0$ with index $\alpha$,   {\it i.e.}
\[
\frac{1-f(x^*-tx)}{1-f(x^*-t)} \rightarrow x^{\alpha}, \mbox{ as } t \rightarrow 0,
\]
for some $\alpha \in \bR$, then for every $\epsilon>0$, there exists $t_0>0$ such that,  $\forall 0< t < t_0$, $\forall 0< x <1$,
\begin{equation} \label{Potterf}
(1-\epsilon) x^{\alpha+\epsilon} \leq \frac{1-f(x^*-tx)}{1-f(x^*-t)} \leq (1+\epsilon) x^{\alpha-\epsilon}
\end{equation}
and
\begin{equation} \label{Potter1surf}
(1-\epsilon) x^{-\alpha+\epsilon} \leq \frac{1-f(x^*-t)}{1-f(x^*-tx)} \leq (1+\epsilon) x^{-\alpha-\epsilon}
\vspace{0.4cm}
\end{equation}
\end{cor}

\noindent Below, $U$ corresponds to the quantile function associated to $H$ introduced in paragraph 2.1.

 \begin{cor}
If $U$ satisfies  condition $(\ref{CondU})$, then for every $\epsilon>0$, there exists $t_0>0$ such that,  $\forall 0< t < t_0$, $\forall x \geq 1$,
\begin{equation}  \label{BornesPotterU}
(1-\epsilon) x^{\gamma-\epsilon} \leq \frac{U(\infty)- U(tx)}{U(\infty)- U(t)} \leq (1+\epsilon) x^{\gamma+\epsilon} .
\end{equation}
\end{cor}

 \begin{prop} (see \cite{deHaanFerreira2006}  Theorem B.2.18)
 \label{propLogU}\\
If $U$ satisfies  condition $(\ref{CondLogU})$ with the positive function $a$, then there exists a  function $q_0$ equivalent to $a/U$ at infinity such that
$\forall  \epsilon >0$, $\exists t_0 >0$, $\forall t \geq  t_0$, $\forall x \geq 1$,
\begin{equation} \label{BornesPotterLogU}
\frac{x^{\gamma}-1}{\gamma} - \epsilon x^{\gamma+\epsilon} \leq \frac{\log U(tx) - \log U(t)}{q_0(t)} \leq \frac{x^{\gamma}-1}{\gamma} + \epsilon x^{\gamma+\epsilon} .
\end{equation}
\end{prop}

\zdeux

We now proceed to the proofs of the different lemmas stated previously and finally of Proposition \ref{lemmeegaliteMnMntilde} .

\subsection{Proof of Lemma \ref{lemmeXin}}

Let $c_n := \left( q_0(n/k_n)/a_{n,k} \right)^{\alpha}$ (which tends to $1$ as $n\ra \infty$) and  $(Y_i)$ be a sequence of i.i.d standard Pareto random variables. Let
 \[
 LL_{i,k} :=  \frac{\log(U(Y_{n-i+1,n})) -\log(U(Y_{n-k_n,n}))}{q_0(Y_{n-k_n,n})}
  \mbox{ and  } QQ_{i,k} := \frac{\left( \frac{Y_{n-i+1,n}}{Y_{n-k_n,n} }\right)^\gamma -1}{\gamma} .
  \]
For every $1\leq i\leq k_n$, we thus have
\[
 \begar{lll}
 \displaystyle c_n^{-1}\frac{\xi_{i,n}}{a_{n,k}^{\alpha}}
 & \stackrel{d}{=}  &  i \left( (LL_{i,k})^{\alpha} - (LL_{i+1,k})^{\alpha} \right) \zun \\
 & =  &  \alpha_{i,n} + \beta_{i,n},
 \finar  \vspace{-0.2cm}
\]
where \vspace{-0.2cm}
\begin{eqnarray*}
  \alpha_{i,n} & := & i \left( (QQ_{i,k})^{\alpha}- (QQ_{i+1,k})^{\alpha} \right) \zun\\
  \beta_{i,n}  & := & i (B_{k,n}(i) - B_{k,n}(i+1)) \zun\\
  B_{k,n}(i) & := & (LL_{i,k})^{\alpha} - (QQ_{i,k})^{\alpha} ,
\end{eqnarray*}
with $B_{k,n}(k_n+1)=0$.
\ztrois

Our first step will be to prove that (with $\delta'>0$)
\[
\mbox{  \ } \frac 1 {k_n} \sum_{i=1}^{k_n} u_i^{-a} \beta_{i,n}
\; \mbox{ or } \;
\frac{1}{k_n^{a+\delta'}} \sum_{i=1}^{k_n}  u_i^{-a} \beta_{i,n}
\ \mbox{ is $ o_{\bP}(1)$ whether $0<a<1$ or $a > 1$}.
\]
Using bounds $(\ref{BornesPotterLogU})$ for some $\epsilon' >0$, with $t=Y_{n-k_n,n}$ and $x= Y_{n-i+1,n}/Y_{n-k_n,n}  >1$, and relying on the mean value theorem, we easily prove that, since $\alpha \geq 1$ and $\gamma<0$,
\begeq{LL-QQ}
|(LL_{i,k})^{\alpha}- (QQ_{i,k})^{\alpha}| \ \leq \ c \; \epsilon' \ \left(\frac{Y_{n-i+1,n}}{Y_{n-k_n,n}} \right)^{\gamma+ \epsilon'} \ \leq \ c\;\epsilon'
\fineq
 for some constant $c$ (close to $\alpha |\gamma|^{1-\alpha}$).  Therefore  $|B_{k,n}(i)| = o_{\bP}(1)$,  uniformly on $1 \leq i \leq k_n$.  Since $|B_{k,n}(i)| = \left| \sum_{j=i}^{k_n}  \frac{\beta_{j,n}}{j} \right|$, we thus have, when $0<a<1$,
 \begin{eqnarray*}
\left| \textstyle \frac 1 {k_n+1} \sum_{i=1}^{k_n} u_i^{-a} \beta_{i,n} \right|
& = &
 \textstyle \left| \sum_{i=1}^{k_n}  \frac{\beta_{i,n}}{i} u_i^{1-a}\, ds \right|\\
& = &
 \textstyle(1-a) \left| \sum_{i=1}^{k_n}  \frac{\beta_{i,n}}{i} \int_0^{u_i} s^{-a}\, ds \right|\\
& \leq &\textstyle (1-a) \sum_{i=1}^{k_n} \left| \sum_{j=i}^{k_n}  \frac{\beta_{j,n}}{j} \right| \int_{u_{i-1}}^{u_i} s^{-a}\, ds \\
& = &
\textstyle o_{\bP}(1)  \sum_{i=1}^{k_n}  (u_i^{1-a} - u_{i-1}^{1-a}) \\
& \leq &
\textstyle o_{\bP}(1)  \frac 1{k_n} \sum_{i=1}^{k_n}  u_i^{-a} \\
& = & o_{\bP}(1) .
\end{eqnarray*}
The proof for  $a >1$ is similar (see end of subsection 5.2.1 in \cite{WWExtremes2014} for more details, with the difference that now $|B_{k,n}(i)|=o_{\bP}(1)$ uniformly in $i$).
\zcinq \\
Since we have dealt with the $\beta_{i,n}$ part, the lemma will be proved as soon as we obtain that, when $0<a<1$,
\begeq{convAlphain_1}
 \frac 1{k_n} \sum_{i=1}^{k_n}  u_i^{-a}  \alpha_{i,n}  \ \limprob \ (1-a) \  |\gamma|^{-\alpha-1} Beta\left( \, \frac{1-a}{|\gamma|} \, ; \,  \alpha+1 \, \right)
\fineq
and, when $a >1$,
\begeq{convAlphain_2}
 \frac{1}{k_n^{a+\delta'}} \sum_{i=1}^{k_n}  u_i^{-a} \alpha_{i,n}  \ \limprob \  0.
\fineq
From now on we will sometimes  write $k$ instead of $k_n$. Let $(E_{i})$ be a sequence of i.i.d standard exponential random variables. According to $(\ref{loiY})$, and by applying the mean value theorem, there exist some random variables $E^*_{i,k}\in [E_{k-i,k},E_{k-i+1,k}]$ such that (remind below that $\gamma$ is $<0$ and $\alpha\geq 1$)
\begin{eqnarray}
 \alpha_{i,n}
  &  \stackrel{d}{=} &  \textstyle
  i \left( \left(  \frac{e^{\gamma E_{k-i+1,k}}-1}{\gamma} \right)^{\alpha}  -  \left(  \frac{e^{\gamma  E_{k-i,k}}-1}{\gamma} \right)^{\alpha} \right)
 \zun \nonumber \\
  & = & \textstyle
  \alpha \ i \, (E_{k-i+1,k} - E_{k-i,k}) \  e^{\gamma  E^*_{i,k}} \left(  \frac{1-e^{\gamma E^*_{i,k}}}{|\gamma|} \right)^{\alpha-1}
 \zun \nonumber \\
 & = & \alpha|\gamma|^{1-\alpha} \times i \, (E_{k-i+1,k} - E_{k-i,k}) \left\{ \ u_i^{|\gamma|}(1-u_i^{|\gamma|})^{\alpha-1} \ + \ \Delta_{i,n} \ \right\} \label{decompAlphain}
\end{eqnarray}
where
\[
 \Delta_{i,n} :=e^{\gamma  E^*_{i,k}} (1-e^{\gamma E^*_{i,k}})^{\alpha-1} - u_i^{|\gamma|} (1- u_i^{|\gamma|})^{\alpha-1} .
\]
We will prove later that
\begeq{convDeltain}
\max_{i\leq k_n} |\Delta_{i,n}|=o_{\bP}(1).
\zun
\fineq
For the moment, note that  $ (\, i \ (E_{k-i+1,k} - E_{k-i,k}) \, )_{i\leq k_n} \stackrel{d}{=} (f_i)_{i\leq k_n}$ due to the Renyi representation, where $(f_i)$ denotes a sequence of i.i.d standard exponential random variables. Moreover, application of the law of large numbers for triangular arrays of independent random variables  (cf \cite{ChowTeicher1997} ; details are omitted) implies that, when $0<a<1$, \small
\begeq{chowteicher1}
\frac 1 {k_n} \sum_{i=1}^{k_n} u_i^{-a}  f_i= O_{\bP}(1)   \makebox[1.2cm][c]{ and }
\frac 1{k_n} \sum_{i=1}^{k_n}  u_i^{-a} f_i \; u_i^{|\gamma|}(1-u_i^{|\gamma|})^{\alpha-1}  \limprob \int_0^1 \ x^{|\gamma|-a} (1- x^{|\gamma|})^{\alpha-1} \ dx
\fineq
\normalsize and,  when $a>1$ (and $\delta'>0$ is given),
\begeq{chowteicher2}
 \frac 1 {k_n^{a+\delta'}} \sum_{i=1}^{k_n} u_i^{-a}  f_i= o_{\bP}(1) .
\fineq
Considering first the situation $0<a<1$, combining (\ref{decompAlphain}), (\ref{convDeltain}) and  (\ref{chowteicher1}) shows that relation (\ref{convAlphain_1}) will hold as soon as
\[
  \int_0^1 \ x^{|\gamma|-a} (1- x^{|\gamma|})^{\alpha-1} \ dx = \alpha^{-1} \,(1-a) \,  |\gamma|^{-2} \, Beta\left( \, \frac{1-a}{|\gamma|} \, ; \,  \alpha+1 \, \right) .
\]
Use of the formulas $Beta(u,v)=\Gamma(u)\Gamma(v)/\Gamma(u+v)$ and $u\Gamma(u)=\Gamma(u+1)$ proves the latter relation.
When $a>1$,  relation (\ref{convAlphain_2}) is a consequence of (\ref{decompAlphain}), (\ref{convDeltain}) and  (\ref{chowteicher2}) .
\bigskip

It remains to prove relation (\ref{convDeltain}) . For this purpose, we introduce the sequence $(V_i)$ of i.i.d. standard uniform random variables such that $V_{i,k} = e^{-E_{k-i+1,k}}$, and we note $V^*_{i,k}:= e^{-E^*_{i,k}}$. If we set $W_{i} := V^{|\gamma|}_{i,k}$, $W^*_{i} := (V^*_{i,k})^{|\gamma|}$ and $v_i := u_i^{|\gamma|}$, then relation (\ref{convDeltain}) is now
 \[
 M_n := \max_{1 \leq i \leq k_n} \left| W^*_{i} (1- W^*_{i})^{\alpha-1}- v_i (1- v_i)^{\alpha-1}  \right| = o_{\bP} (1),
 \]
 with $0 < W_{i} \leq W^*_{i} \leq W_{i+1}<1$. Unfortunately, the function $x\mapsto x(1-x)^{\alpha-1}$ is not uniformly continuous on $]0,1[$ if $\alpha$ is smaller than 2. Until the end of the proof we will note $\max_i$ instead of $\max_{1\leq i\leq k_n}$. Since $|x-a|\leq \max\{|y-a|,|z-a|\}$ whenever $|x-y|\leq |z-y|$ (this yields the first inequality below), we have \small
\[
 \begar{ll}
 M_n  &
 \leq  \max \left\{ \,
 \max_i | W_{i} (1- W_{i+1})^{\alpha-1}- v_i (1- v_i)^{\alpha-1} |
 \ ; \  \right.
 \\ &
    \hspace*{4.3cm}\left.
 \max_i | W_{i+1} (1- W_{i})^{\alpha-1}- v_i (1- v_i)^{\alpha-1} |
 \, \right\}
 \zdeux\\ &
 \leq  \max \left\{ \,
 \max_i | W_{i} (1- W_{i+1})^{\alpha-1}- W_i (1- W_i)^{\alpha-1} |
 \ ; \  \right.
 \\ &
    \hspace*{1.3cm}\left.
 \max_i | (W_{i+1}-W_i) (1- W_i)^{\alpha-1}  |
 \ ; \;
 \max_i | W_{i} (1- W_{i})^{\alpha-1}- v_i (1- v_i)^{\alpha-1} |
 \, \right\}
 \zdeux\\ &
 \leq   \max \left\{ \,
 \max_i | (1- W_{i+1})^{\alpha-1} - (1- W_i)^{\alpha-1} |
 \ ; \  \max_i | W_{i+1}-W_i |
 \ ; \ \right.
 \\ &
    \hspace*{2.3cm}\left.
 \max_i | W_{i}-v_i |
 \ ; \;
 \max_i | (1- W_{i})^{\alpha-1} - (1- v_i)^{\alpha-1} |
 \, \right\}
 \zdeux\\ &
 =: \max \left\{ \ M_{n,1} \ ; \ M_{n,2} \ ; \ M_{n,3} \ ; \ M_{n,4} \   \right\} \, .
 \finar
\]
\normalsize Now, since $|x^{a}-y^{a}| \leq (a\wedge 1)|x-y|^{a\vee 1}$ whenever $x$ and $y$ belong to $[0,1]$ and $a>0$, we have $M_{n,1}\leq cM_{n,2}^{c'}$ and $M_{n,4}\leq cM_{n,3}^{c'}$ for some positive constants $c$ and $c'$ . On the other hand, $M_{n,2}$ is bounded by $2M_{n,3}+\max_i|v_{i+1}-v_i| = 2M_{n,3} + o(1)$. Therefore, the negligibility of $M_n$ amounts to
$$
M_{n,3}=\max_{1\leq i\leq k_n} \left|V_i^{|\gamma|} - u_i^{|\gamma|}\right| = o_{\bP}(1).
$$
This property is proved in details in \cite{Beirlant2002} (page 164, with $-\rho$ instead of $|\gamma|$), so we do not reproduce it here.
\hfill$\diamond$

 \subsection{Proof of Lemma \ref{lemmeBin}}

If $Y_{1,n}, \ldots, Y_{n,n}$ denote the ascending order statistics of $n$ i.i.d standard Pareto random variables, we have
\[
   \tilde{Z}_{i,n} \stackrel{d}{=} \frac{U(\infty)-U(Y_{n-i+1,n})}{U(\infty)-U(Y_{n-k_n,n})} .
\]
Applying bounds $(\ref{BornesPotterU})$, it comes,  for some given $\epsilon'>0$  and $n$ sufficiently large,
\begin{eqnarray*}
\left( \frac{U(\infty)-U(Y_{n-i+1,n})}{U(\infty)-U(Y_{n-k_n,n})} \right)^{\theta}
 & \leq &  (1+\epsilon')^{\theta} \left(\frac{Y_{n-i+1,n}}{Y_{n-k_n,n}}\right)^{a_+} \\
 \left( \frac{U(\infty)-U(Y_{n-i+1,n})}{U(\infty)-U(Y_{n-k_n,n})} \right)^{\theta'}
 & \geq &  (1-\epsilon')^{\theta'} \left(\frac{Y_{n-i+1,n}}{Y_{n-k_n,n}}\right)^{a_-}
\end{eqnarray*}
where
$$
 a_+ = (\gamma+\epsilon')\theta \ \mbox{\ and \ } \ a_-=(\gamma-\epsilon')\theta'.
 $$
We finish  the proof as for  Lemma 1 of \cite{WWExtremes2014}.
\hfill$\diamond$

\subsection{Proof of Proposition \ref{lemmeegaliteMnMntilde}}
(i) First, note that for any sequences $(a_n)$ and $(b_n)$ such that  $a_0=0$ and $b_{k+1}=0$, we have $\sum_{i=1}^k a_i (b_i-b_{i+1})= \sum_{i=1}^k  (a_i - a_{i-1}) b_i$.
By letting $a_i=\frac{i}{1-\hat G_n(Z^-_{n-i+1,n})} $ and $b_i=\log^{\alpha} \left(  \frac{Z_{n-i+1,n}}{Z_{n-k_n,n}} \right)$, we have $b_{k+1}=0$  but $a_0$ is undefined, so we set $a_0=0$. Therefore, by the definition (\ref{defxiin}) of $\xi_{i,n}$, this implies that
\[
\begar{rcr}
\sum_{i=1}^{k_n} \frac{\xi_{i,n}}{1-\hat G_n(Z^-_{n-i+1,n})}
& = & \sum_{i=2}^{k_n} \left( \frac{i}{1-\hat G_n(Z^-_{n-i+1,n})}  -  \frac{i-1}{1-\hat G_n(Z^-_{n-i+2,n})}  \right) \   \log^{\alpha} \left(  \frac{Z_{n-i+1,n}}{Z_{n-k_n,n}} \right)    \zun \\ & & +  \frac{1}{1-\hat G_n(Z^-_{n,n})}  \   \log^{\alpha} \left(  \frac{Z_{n,n}}{Z_{n-k_n,n}} \right).
\finar
\]
By definition of $M^{(\alpha)}_{n,k_n}$ and $D^{(\alpha)}_{n,k_n}$, part (i) of Proposition \ref{lemmeegaliteMnMntilde} will be proved as soon as we show that, for $2 \leq i \leq k_n$,
\[
\frac{i}{1-\hat G_n(Z^-_{n-i+1,n})}  -  \frac{i-1}{1-\hat G_n(Z^-_{n-i+2,n})}  = \frac{\delta_{n-i+1,n}}{1-\hat G_n(Z^-_{n-i+1,n})} .
\]
Note that
\begin{eqnarray}
 & \frac{i}{1-\hat G_n(Z^-_{n-i+1,n})}  -  \frac{i-1}{1-\hat G_n(Z^-_{n-i+2,n})}  \nonumber \zun \\
= &  \frac{1}{1-\hat G_n(Z^-_{n-i+1,n})} \left( i - (i-1) \frac{1-\hat G_n(Z^-_{n-i+1,n})}{1-\hat G_n(Z^-_{n-i+2,n})} \right) \label{eg2}
\end{eqnarray}
where
\[
\frac{1-\hat G_n(Z^-_{n-i+1,n})}{1-\hat G_n(Z^-_{n-i+2,n})} = \frac{\prod_{j=1}^{n-i} \left(  \frac{n-j}{n-j+1} \right)^{1-\delta_{j,n}}}{\prod_{j=1}^{n-i+1} \left(  \frac{n-j}{n-j+1} \right)^{1-\delta_{j,n}}}.
\]
The right-hand side is clearly equal to $\frac{i}{i-1} \left( \frac{i-1}{i} \right)^{\delta_{n-i+1,n}}$, which may be rewritten as $\frac{i}{i-1}  \left( 1-\frac{\delta_{n-i+1,n}}{i} \right)$. Replacing in  $(\ref{eg2} )$ concludes the proof of (i).
\zdeux

\noindent (ii) Note first that
 \[
 D^{(\alpha)}_{n,k_n} = \frac{1}{k_n} \ A_n \ C_{1,n}  W_{1,n},
 \]
 where $A_n$ and $C_{i,n}$ are defined in the beginning of the proof of Theorem \ref{consistKSV} and
 \[
 W_{1,n} = \frac{1- G(Z_{n-k_n,n})}{1- G(Z_{n,n})}  \log^{\alpha} \left( \frac{Z_{n,n}}{Z_{n-k_n,n}} \right).
 \]
 We know that $A_n=o_{\bP}(1)$. Since $C_{1,n}=O_{\bP}(1)$ (see Theorem 2.2 in  \cite{Zhou1991}), we only have to prove that $ \frac{1}{k_n}   \frac{W_{1,n}}{a^{\alpha}_{n,k}}=o_{\bP}(1)$.
For the same reasons that led to  $(\ref{bornesWin})$, for any $\epsilon >0$ and $n$ sufficiently large,
  we get $W_{1,n} \leq  W^+_{1,n}$ with
\[
W^+_{1,n} = (1+\epsilon)  \log^{\alpha} \left( \frac{Z_{n,n}}{Z_{n-k_n,n}} \right) \tilde{Z}_{1,n}^{\gamma_C^{-1} - \epsilon}.
\]
If $(Y_i)_{1\leq i \leq n}$ is an i.i.d. sequence of standard Pareto random variables, then
\[
c_n^{-1} \frac{W^+_{1,n}}{a^{\alpha}_{n,k}} \stackrel{d}{=} (1+\epsilon) (LL_{1,k})^{\alpha} (U_{1,k})^{\gamma_C^{-1} - \epsilon},
\]
where $c_n \rightarrow 1$ and $LL_{1,k}$ were  defined  in the proof of Lemma \ref{lemmeXin} and $U_{1,k}=  \frac{U(\infty)- U(Y_{n,n})}{U(\infty)- (Y_{n-k_n,n})}$.  On one hand,  Potter bounds $(\ref{BornesPotterU})$ yield for any $\epsilon' >0$ and $n$ sufficiently large,
\[
(U_{1,k})^{\gamma_C^{-1} - \epsilon} \leq (1+\epsilon') \left( \frac{Y_{n,n}}{Y_{n-k_n,n}}  \right)^{a^+},
\]
with $a^+=(\gamma+\epsilon')(\gamma_C^{-1} - \epsilon)$.
On the other hand, $(\ref{LL-QQ})$ implies that  for $\epsilon''>0$, $|(LL_{1,k})^{\alpha}-(QQ_{1,k})^{\alpha}| \leq c \epsilon''$, for some constant $c$, where $QQ_{i,k}$ is also defined in  the proof of Lemma $\ref{lemmeXin}$. Moreover, it is known that
\begeq{loiY}
 \left(Y_{n-i+1,n} / Y_{n-k_n,n}\right)_{1\leq i\leq k_n}  \stackrel{d}{=} (\tilde Y_{k_n-i+1,k_n})_{1\leq i\leq k_n},
 \fineq
where $\tilde Y_{1,k_n}, \ldots, \tilde Y_{k_n,k_n}$ are the ascending order statistics of $k_n$ i.i.d random variables $\tilde Y_1,\ldots,\tilde Y_{k_n}$
with standard Pareto distribution. \zun

Consequently,
$(LL_{1,k})^{\alpha} (U_{1,k})^{\gamma_C^{-1} - \epsilon} \leq Q_{1,n}$,
with
\[
Q_{1,n} \stackrel{d}{=}  (1+\epsilon') \left(  \frac{\tilde Y_{k_n,k_n}-1}{\gamma}\right)^{\alpha}  (\tilde Y_{k_n,k_n})^{a^+}  + c \epsilon''  (1+\epsilon')  (\tilde Y_{k_n,k_n})^{a^+} .
\]
Since $\tilde Y_{k_n,k_n} \geq 1$  and $\gamma <0$, the right-hand side of the relation above  is lower  than $c' (\tilde Y_{k_n,k_n})^{a^+} $, for some constant $c'$. \zun

Now, Standard Pareto distributions having moments of order less than $1$,
\[
\tilde Y_{k_n,k_n} = \max_{1 \leq i \leq k_n} \tilde Y_{i,k_n} = o_{\bP}(k_n^p),
\]
for any $p >1$, therefore, $\frac{1}{k_n} (\tilde Y_{k_n,k_n})^{a^+} = o_{\bP}(k_n^{pa^+-1})$, which is clearly  $o_{\bP}(1)$. This concludes the proof.
\hfill$\diamond$


\vspace{1.cm}

\begin{thebibliography}{15}

\bibitem[Beirlant et~al.~(2002)]{Beirlant2002}
 J. Beirlant, G. Dierckx, A. Guillou and C. St\v{a}ric\v{a} .
\newblock On exponential representations of log spacings of order statistics.
\newblock  In {\em Extremes} {\bf 5}, pages 157-180 (2002)

\bibitem[Beirlant et~al.~(2007)]{Beirlant2007}
 J. Beirlant, G. Dierckx, A. Guillou and A. Fils-Villetard .
\newblock Estimation of the extreme value index and extreme quantiles under random censoring.
\newblock  In {\em Extremes} {\bf 10}, pages 151-174  (2007)

\bibitem[Beirlant et~al.~(2010)]{Beirlant2010}
 J. Beirlant, A. Guillou and G. Toulemonde .
\newblock Peaks-Over-Threshold modeling under random censoring.
\newblock  In {\em Comm. Stat. : Theory and Methods} {\bf 39}, pages 1158-1179 (2010)

\bibitem[Brahimi et~al.~(2013)]{BrahimiMeraghniNecir2013}
 B. Brahimi, D. Meraghni and A. Necir .
\newblock Approximations to the tail index estimator of a heavy-tailed distribution under random censoring and application.
\newblock  working paper, arXiv:1302.1666  (2013)

\bibitem[Chow  and Teicher (1997)]{ChowTeicher1997}
Y.S. Chow and H. Teicher .
\newblock Probability theory. Independence, interchangeability, martingales.
\newblock {\em Springer}   (1997)

\bibitem[Cs\"org\H{o} (1996)]{Csorgo96}
S. Cs\"org\H{o} .
\newblock Universal Gaussian approximations under random censorship.
\newblock In {\em Annals of statistics} {\bf 24 (6)}, pages 2744-2778 (1996)

\bibitem[Diop et~al.~(2014)]{DupuyDiopNdao2014}
A. Diop, J-F. Dupuy and P. Ndao.
\newblock Nonparametric estimation of the conditional tail index and extreme quantiles under random censoring.
\newblock In {\em Computational Statistics \& Data Analysis} {\bf 79}, pages 63-79 (2014)



\bibitem[Einmahl et~al.~(2008)]{Einmahl2008}
J. Einmahl, A. Fils-Villetard   and  A. Guillou .
\newblock Statistics of extremes under random censoring.
\newblock In {\em Bernoulli} {\bf 14}, pages 207-227 (2008)


\bibitem[Gomes and Neves (2011)]{GomesNeves2011}
M.I. Gomes and M.M. Neves .
\newblock Estimation of the extreme value index for randomly censored data.
\newblock In {\em Biometrical Letters} {\bf 48 (1)}, pages 1-22 (2011)

\bibitem[Haan and Ferreira (2006)]{deHaanFerreira2006}
L. de Haan and A. Ferreira  .
\newblock Extreme Value Theory : an Introduction.
\newblock {\em Springer Science + Business Media} (2006)

\bibitem[Segers (2001)]{Segers2001}
J. Segers.
\newblock Residual estimators.
\newblock In {\em Journal of Statistical Planning and Inference} {\bf 98}, pages 15-27 (2001)







\bibitem[Worms (2014)]{WWExtremes2014}
J. Worms and R. Worms
\newblock New estimators of the extreme value index under random right censoring, for heavy-tailed distributions.
\newblock In {\em Extremes}  {\bf 17 (2)}, pages 337-358  (2014)


\bibitem[Zhou (1991)]{Zhou1991}
M. Zhou.
\newblock Some Properties of the Kaplan-Meier Estimator for Independent Nonidentically Distributed Random Variables.
\newblock In {\em Annals of statistics} {\bf 19 (4)}, pages 2266-2274 (1991)





\end{thebibliography}
\end{document}